\documentclass[11pt]{article}
\usepackage{amssymb,color}
\usepackage{amsmath}
\usepackage[pctex32]{graphics}
\usepackage{graphicx}

  \textwidth
14.5 cm \textheight 22.2 cm 

\def\be{\begin{equation}}
\def\ee{\end{equation}}
\newcommand{\fer}[1]{\eqref{#1}}
\def\dis{\displaystyle}

\def\R{\rho}                          
\def\A{\mathfrak{A}}

\newcommand{\p}{\partial}                 
\def\N{\mathbb{N}}                        
\def\ff{\hat f}                           
\def\Realr{\mathbb R}                        
\def\rr{\mathbb R}                        
\newcommand{\Rn}{{{\rr}^N}}               
\newcommand{\Rno}{\Rn\hspace{-1.5mm}_o}   
\newcommand{\Rto}{\rr^3\hspace{-1.5mm}_o} 
\def\et{\tilde{e}}                    
\def\tQ{\tilde{Q}}                    

\renewcommand{\P}{{\cal P}}                
\newcommand{\Pn}{\P(\Rn)}                  
\newcommand{\Ptwo}{\P_2(\Rn)}              
\newcommand{\Ptwot}{\P_2(\rr^3)}           
\def\expec{\mathbb E}                      
\newcommand{\Rop}{{\mathbf R}}
\newcommand{\setR}{{\mathbb R}}
\newcommand{\setC}{{\mathbb C}}

\newtheorem{theorem}{Theorem}[section]
\newtheorem{lemma}[theorem]{Lemma}
\newtheorem{proposition}[theorem]{Proposition}
\newtheorem{corollary}[theorem]{Corollary}
\newtheorem{remark}[theorem]{Remark}

\def\proof{\par{\it Proof.-} \ignorespaces}
\def\endproof{{\ \vbox{\hrule\hbox{%
   \vrule height1.3ex\hskip0.8ex\vrule}\hrule
  }}\par}

\begin{document}

\title{Over-populated Tails for conservative-in-the-mean \\
Inelastic Maxwell Models}

\author{ J. A. Carrillo
\thanks{ICREA-Departament de Matem\`atiques,
  Universitat Aut\`onoma de Barcelona,
  E-08193 Bellaterra, Spain.
  E-mail: {\tt carrillo@mat.uab.es}.} \,
S. Cordier
\thanks{F\'ed\'eration Denis Poisson (FR 2964),
Department of Mathematics (MAPMO UMR 6628)
University of Orl\'eans and CNRS,
F-45067 Orl\'eans, France.
E-mail:{\tt Stephane.Cordier@univ-orleans.fr} }\,
G. Toscani
\thanks{Dipartimento di Matematica,
Universit\`a di Pavia, via Ferrata 1,
I-27100 Pavia, Italy.
E-mail: {\tt giuseppe.toscani@unipv.it}.}}

\maketitle


\newcommand{\which}{that }

\begin{abstract}
We introduce and discuss spatially homogeneous Maxwell-type models of the
nonlinear Boltzmann equation undergoing binary collisions with a random
component. The random contribution to collisions is such that the usual
collisional invariants of mass, momentum and energy do not hold pointwise, even
if they all hold in the mean. Under this assumption it is shown that, while the
Boltzmann equation has the usual conserved quantities, it possesses a steady
state with power-like tails for certain random variables. A similar situation
occurs in kinetic models of economy recently considered by two of the authors
\cite{MaTo06}, which are conservative in the mean but possess a steady
distribution with Pareto tails. The convolution-like gain operator is
subsequently shown to have good contraction/expansion properties with respect
to different metrics in the set of probability measures. Existence and
regularity of isotropic stationary states is shown directly by constructing
converging iteration sequences as done in \cite{Bobylev-Cercignani2}.
Uniqueness, asymptotic stability and estimates of overpopulated high energy
tails of the steady profile are derived from the basic property of
contraction/expansion of metrics. For general initial conditions the solutions
of the Boltzmann equation are then proved to converge with computable rate as
$t\to\infty$ to the steady solution in these distances, which metricizes the
weak convergence of measures. These results show that power-like tails in
Maxwell models are obtained when the point-wise conservation of momentum
and/or energy holds only globally.
\end{abstract}


\section{Introduction}

In this paper, we introduce and discuss the possibility to obtain
steady solutions with power-like tails starting from conservative
molecular systems described by the Boltzmann equation with
Maxwell-type collision kernels. The starting point of our model is
to consider binary collisions \which result in a linear combination
of an inelastic collision and a random contribution. As we
shall see, the random addition to the post-collision velocities
can only increase the mean of the collisional energy, and, among
other things, it gives the possibility to construct a binary collision
\which preserves (in the mean) mass, momentum and energy. Our model
is closely related to a kinetic model for economics introduced by 
Pareschi and two of the present authors \cite{CPT}. There, the
random contribution to the collision (trade) was introduced to
take into account the returns of the market.

Inelastic Maxwell models were introduced by Bobylev, Gamba and one
of the authors in 2000 \cite{Bobylev-Carrillo-Gamba}; see also
\cite{BK00} for the one dimensional case. Maybe the most
interesting result (absent in the elastic case) is the existence
of self-similar solutions in the homogeneous cooling problem and
the non-Maxwellian behavior of these solutions, which displays
power-like decay for large velocities. It was conjectured in
\cite{EB2} and later proved in \cite{Bobylev-Cercignani2,
Bobylev-Cercignani-Toscani, BisiCT2} that such solutions represent
intermediate asymptotics for a wide class of initial data. Other
results concerned with self-similar solutions in the theory of
the classical (elastic) Boltzmann equation for Maxwell molecules were
also recently published in \cite{Bobylev-Cercignani,
Bobylev-Cercignani2}.  In light of these results, it looks
clear that in many aspects both elastic and inelastic Maxwell
models must be studied from a unified point of view. As observed
by Bobylev and Gamba in
\cite{Bobylev-Cercignani-Gamba,Bobylev-Cercignani-Gamba2}, an
interesting question arises in connection with power-like tails
for high velocities. Is it possible to observe a similar effect,
i.e., an appearance of power-like tails from initial data with
exponential tails, in a system of particles interacting according
to laws of classical mechanics without energy loss? In
\cite{Bobylev-Gamba} Bobylev and Gamba gave a partial answer to
this question by showing that, under a certain limiting procedure,
such behavior can in principle be observed if one considers a
mixture of classical Maxwell gases. More precisely, self-similar
solutions converging towards maxwellian equilibrium were proved to
have power-like tails once normalized by the equilibrium.

In this paper, we will try to elucidate the same question,
starting from a somewhat different point of view. Our starting
point will be a suitable modification to the homogeneous Boltzmann
equation for the inelastic Maxwell molecules introduced in Ref.
\cite{Bobylev-Carrillo-Gamba}, in such a way that the usual
conservations of mass, momentum and energy in the binary
collisions still continue to hold in the mean sense. The
scaled-in-time inelastic Boltzmann equation introduced in
\cite{Bobylev-Carrillo-Gamba} reads
 \begin{equation}
\frac{\p f}{\p t}=Q_e(f,f)\,. \label{BE1}
 \end{equation}
Here, $f(v,t)$ is the density for the velocity space distribution
of the molecules at time $t$, while $Q_e(f,f)$ is the inelastic
Boltzmann collision operator, which contains the effects of binary
collisions of grains. As usual in this context, the collision
operator $Q_e(f,f)$ is more easily treated if expressed in weak
form. This corresponds to writing, for every suitable test
function~$\varphi$,
  \begin{equation}
(\varphi,Q_e(f,f))= \frac{1}{4 \pi} \int_{\Realr^3}
\int_{\Realr^3} \int_{S^2} f(v) f(w) \Big[ \varphi (v^*)- \varphi
(v) \Big] dv\,dw\,d\sigma. \label{Qweak}
 \end{equation}
In \fer{Qweak}, $v^*$ is the outgoing velocity corresponding to a
particle in the collision defined by the incoming velocities $v,w$
and the angular parameter $\sigma \in S^2$:
\be
 \begin{split}
 &\dis v^*= \frac{1}{2} (v+w)+ \frac{1-e}{4}(v-w)+\frac{1+e}{4}
|v-w|\sigma,
\\*[.3cm]
&\dis w^*= \frac{1}{2} (v+w)- \frac{1-e}{4}(v-w)-\frac{1+e}{4}
|v-w|\sigma\,.
\end{split}
\label{post}
 \ee
 The parameter $0\leq e \leq 1$ represents the restitution coefficient.

In the model we consider, this restitution coefficient will be chosen as a
random variable \which can be interpreted from a physical point of view as the
stochasticity in the microscopic process of collision due to the randomness of
the grains' geometry and the mechanical properties of the medium. We will show in
the next section that this random behavior in restoring energy leads to a precise
form of the energy gain term \which differs from the usually chosen diffusion
term, the so-called "thermal bath". This new form of "thermal bath" is thus related
to the process generated by the randomness of the granular media. We prove
that this particular thermal bath yields equilibrium states with power law
tails.

Such over-populated tails in distributions at equilibrium arise in
other contexts.   We shall present similar results on the large
time behavior of collisional kinetic theory applied to economic
modelling. In this framework, the kinetic variable represents the
wealth of agents and the collision operator describes the
evolution of the wealth distribution through exchanges. We refer
to \cite{CPT, MaTo06} and references therein for a mathematical
presentation of these models closely related to so called
"econo-physics". In such models, the equations between pre- and
post-collisional values involve some randomness \which is related
to the stochasticity of the market that provides random returns.


In the remainder of this paper, we will study in detail the large time
behavior of the solution of the Boltzmann equation involving such a stochastic
process. We show that the validity (at a macroscopic level) of the classical
collision invariants is enough to guarantee convergence towards a steady
profile, but not enough to reach a Maxwellian-like profile. In fact, we will
show that there is a class of random perturbations of the coefficient of
restitution such that the steady state possesses power-like tails.

A crucial role in our analysis is played by the weak norm
convergence, which is obtained by further pushing the development
of a method first used in \cite{Gabetta-Toscani-Wennberg} to
control the exponential convergence of Maxwellian molecules in
certain weak norms. This will be done by using the fact that the
nonlinear operator in the Boltzmann equation (see \fer{BE2})
 can be expressed in Fourier
variables in a simple closed form using Bobylev's identity
\cite{Bobylev}. Estimates of the evolution of the Wasserstein
distance \cite{Wasserstein69,Vil03,Vil06} between solutions will
be presented for the economic and the inelastic model since they
give complementary information with respect to the results in
\cite{MaTo06}. Concerning this second aspect, we will take
advantage of the recent analysis of Bolley and Carrillo
\cite{Bolley-Carrillo, CTPE} of the inelastic Boltzmann equation
for Maxwell molecules. From this analysis, we will obtain the
uniqueness and asymptotic stability of stationary states for this
model. Finally, the appearance of power-like tails for the
asymptotically stable stationary states will be discussed for both
models, giving explicit examples of random variables producing this
behavior.


The paper is organized as follows: in section \ref{Modeling} we
detail the collisional models for both granular media and economy
applications including random coefficients in the relations
between pre- and post-collisional variables. In section
\ref{quick}, we recall the main properties of probability metrics. In
section \ref{eco2}, we investigate large time behavior of the solution
of the kinetic economy model and section \ref{seccontr} is devoted to
large time behavior of stochastic granular media.

Lastly, let us summarize the two main results of this paper :
first, we give some insight into conditions for a collision
operator to lead to power-law tails (conservatism in mean being
some kind of necessary condition); second,  we propose a new form
for the thermal bath with a physically relevant origin (the
restitution coefficient  taking into account the randomness of
granular media).


\section{Modelling issues and diffusion approximation}
\label{Modeling}


Let us present the proposed stochastic granular model
(with a random restitution coefficient) and its diffusion limit and then
recall briefly the similar analysis for the economy model following
\cite{CPT}.

\subsection{Stochastic granular media}

Considering the weak formulation \fer{Qweak},
easy computations show that $(\varphi(v),Q_e(f,f))=0 $ whenever
$\varphi(v) =1$ and $\varphi(v) =v$, while
$(\varphi(v),Q_e(f,f))<0 $ if $\varphi(v) =v^2$. This corresponds
to conservation of mass and momentum, and, respectively, to loss
of energy for the solution to equation \fer{BE1}. For this reason,
if we fix the initial data to be a centered probability density
function, the solution will remain centered at any subsequent time
$t>0$. The loss of energy in a single collision
with a constant restitution coefficient $e$ is given by
 \be
 \label{loss}
 |v'|^2+ |w'|^2 = |v|^2 +|w|^2 -\frac{1-e^2}4 \left( |v-w|^2 -|v-w|(v-w)\cdot
 \sigma \right).
 \ee
The previous formula is the key to our modification of the collisions.
Let us replace the constant coefficient of restitution $e$ with
a stochastic coefficient of restitution $\et$, such that for a
given random variable $\eta$
 \be\label{re}
 \et = e +\eta, \qquad \mbox{with} \quad\langle \eta\rangle =0 \quad \mbox{and} \quad \langle \eta^2 \rangle =
 \beta^2.
 \ee
In \fer{re} and in the rest of the paper, $\langle \cdot \rangle$ denotes the
mathematical expectation of the real-valued random variable $\eta$, i.e.,
integration against a measure $\mu$. For obvious physical reasons, the random
variable $\eta$ has to be chosen to satisfy $\eta \ge -e$, in order to
guarantee that the (random) coefficient of restitution $\et\geq0$. Using $\et$
instead of $e$ in \fer{post} gives that the momentum is {\it conserved in
average} for a suitable choice of the variance. In fact, since
 \be\label{1new}
 \left\langle |v'|^2+ |w'|^2\right\rangle = |v|^2 +|w|^2 -\frac{1-e^2 -\beta^2}4 \left( |v-w|^2 -|v-w|(v-w)\cdot
 \sigma \right),
 \ee
by choosing the variance $\beta^2 = 1-e^2 >0$, we obtain
 \be\label{cons}
\left\langle |v'|^2+ |w'|^2\right\rangle = |v|^2 +|w|^2.
 \ee
 We will call a collision process (or equivalently a random cross section)
 satisfying \fer{cons} \emph{conservative
 in the mean}.
 Let us remark that condition \fer{cons} \emph{cannot be satisfied
if $\et$ takes only values less than $1$}, since in that case
$\et^2$ remains also less than 1 and so does its average
$<\et^2>=e^2+\beta^2 <1$.
 The main idea behind
this is that particles can even gain energy in collisions even though
the total energy is conserved in the mean.

 From the physical point of view, this assumption of energy-gain particle
collisions may seem strange. We will show in the sequel that this energy input
can be interpreted as a sort of thermal bath. Particles are immersed in a medium
that produces this random change in the strength of their relative velocity. We
will argue, based on a derivation of a Fokker-Planck approximation, that this
random component in the collision operator can be approximated by a
second-order differential operator whose diffusion matrix depends on the second
moments of the solution $f$ itself and the random variable $\eta$ (see
\cite{CPT, PaTo06} for a similar approach in one dimension).

This idea allows us to consider a new class of Maxwell-type models, 
from now on called \emph{conservative in the mean}, which are
obtained from (post-collision) velocities given by
\be
 \begin{split}
 &\dis v'= \frac{1}{2} (v+w)+ \frac{1-\et}{4}(v-w)+\frac{1+\et}{4}
|v-w|\sigma,
\\*[.3cm]
&\dis w'= \frac{1}{2} (v+w)- \frac{1-\et}{4}(v-w)-\frac{1+\et}{4}
|v-w|\sigma\,.
\end{split}
\label{postcon}
 \ee
where $\et$ is the random coefficient of restitution defined in \fer{re}, and
$\beta^2 = 1-e^2$. The corresponding Boltzmann equation reads
 \be
\frac{\p f}{\p t}= \tilde{Q}_e(f,f) = \left\langle Q_{\et}(f,f)
\right\rangle\, , \label{BE2}
 \end{equation}
and its corresponding weak form is
  \be
(\varphi,\tilde{Q}_e(f,f))= \frac{1}{4 \pi}\left\langle
\int_{\Realr^3} \int_{\Realr^3} \int_{S^2} f(v) f(w) \Big[ \varphi
(v')- \varphi (v) \Big] dv\,dw\,d\sigma \right\rangle.
\label{Qweak1}
 \end{equation}
In view of our choice of the random contribution to the
coefficient of restitution, we now have
$(\varphi(v),\tilde{Q}_e(f,f))=0 $ whenever $\varphi(v) =1,v,
|v|^2$, that is, the classical collision invariants of the elastic
Boltzmann equation.

\subsection{ Formal diffusive asymptotics }
\label{formal}

Before entering into the study of the large-time behavior of the
Boltzmann equation \fer{BE2}, we shall present here some formal
arguments \which hopefully clarify the action of the random
restitution coefficient in the collision mechanism \fer{postcon}.

To this end, following the same method as in \cite{To-Diss}, letting
$(v',w')$ denote the post-collision velocities \fer{postcon} in our
random collision with $(v^*,w^*)$ as post collision velocities
defined by
the classic inelastic collision \fer{post}, we can split the
velocities into their deterministic and random parts
  \be\label{rel2}
 v' = v^* + \eta \Delta(u, \sigma) \, , \quad w' = w^* - \eta \Delta(u,
 \sigma),
 \ee
  where we let $u = v-w$ and
  \[
 \Delta(u, \sigma) = \frac 14 \left( |u|\sigma - u\right).
  \]
Let us consider a Taylor expansion of $\varphi(v')$ around $\varphi(v^*)$
up to second order in $\eta$.
Thanks to \fer{rel2} we get
\begin{eqnarray}\label{exp}
\varphi(v') = \varphi(v^*) + \eta\,\left(\nabla \varphi(v^*) \cdot
\Delta(u, \sigma)\right) + \frac 12 \eta^2
\sum_{i,j}\frac{\partial^2 \varphi(v^*)}{\partial v^*_i\partial
v^*_j}\Delta_i\Delta_j + \dots
\end{eqnarray}
Thus, taking the mean of the expansion \fer{exp}, and using the property
$\langle \eta \rangle= 0$, we get
 \be\label{exp2}
\langle\varphi(v')\rangle = \varphi(v^*) + \frac 12 \beta^2
\sum_{i,j}\frac{\partial^2 \varphi(v^*)}{\partial v^*_i\partial
v^*_j}\Delta_i\Delta_j + \dots.
 \ee
Truncating the expansion \fer{exp2} after the second--order term and
inserting \fer{exp2} into \fer{Qweak1}, we conclude
 \begin{align}\label{weak2}
(\varphi\, , \, \tilde Q_e(f,f)) \, &\simeq (\varphi\, , \,
Q_e(f,f)) + (\varphi, D_e(f,f)) \\
 &=(\varphi\, , \,  Q_e(f,f))+ \frac {\beta^2}{8\pi}
\int_{\Realr^3}\int_{\Realr^3}\int_{S^2}\sum_{i,j}\frac{\partial^2
\varphi(v^*)}{\partial v^*_i\partial v^*_j}\Delta_i\Delta_j
f(v)f(w)dv\, dw\, d\sigma \,\, .\nonumber
 \end{align}
While the first term in \fer{weak2} $Q_e(f,f)$ is the classical inelastic
Boltzmann collision operator, the second term $D_e(f,f)$ needs to be further
analyzed.

Denoting by $(\null^{*}v, \null^{*}w)$ the pre-collision
velocities in the inelastic collision, and taking into account the
fact that the Jacobian of the transformation
$d\null^*v\,d\null^*w$ into $dv\,dw$ for a constant restitution
coefficient is equal to $e^{-1}$, one obtains
 \begin{align}\label{weak3}
(\varphi\, , \,D_e(f,f)) \,\, &= \frac {\beta^2}{8\pi}
\int_{\Realr^3}\int_{\Realr^3}\int_{S^2}\sum_{i,j}\frac{\partial^2
\varphi(v^*)}{\partial v^*_i\partial v^*_j}\Delta_i\Delta_j
f(v)f(w)dv\, dw\, d\sigma\nonumber \\
&=  \frac {\beta^2}{8\pi}
\int_{\Realr^3}\int_{\Realr^3}\int_{S^2}\frac
1e\sum_{i,j}\frac{\partial^2 \varphi(v)}{\partial v_i\partial v_j}
\null^*\Delta_i\null^*\Delta_j f(\null^*v)f(\null^*w)dv\, dw\,
d\sigma\nonumber \\
&= \frac {\beta^2}{8\pi} \int_{\Realr^3}\,\left[
\sum_{i,j}\frac{\partial^2 \varphi(v)}{\partial v_i\partial v_j}
\int_{\Realr^3}\int_{S^2}\frac 1e \null^*\Delta_i\null^*\Delta_j
f(\null^*v)f(\null^*w)\, dw\, d\sigma\right]\, dv
\nonumber \\
&= \int_{\Realr^3} \!\!\varphi(v) \!\left[\frac
{\beta^2}{8\pi}\sum_{i,j}\frac{\partial^2 }{\partial v_i\partial
v_j}\int_{\Realr^3}\int_{S^2}\frac 1e
\null^*\Delta_i\null^*\Delta_j f(\null^*v)f(\null^*w)\, dw\,
d\sigma\right] dv.
 \end{align}
This shows that, at least for small inelasticity,  the random part of the
collision corresponds to a correction given by the nonlinear diffusion operator
$D_e(f,f)(v)$, where
 \be\label{corr}
  D_e(f,f)(v) =
\frac {\beta^2}{8\pi}\sum_{i,j}\frac{\partial^2 }{\partial
v_i\partial v_j}\int_{\Realr^3}\int_{S^2}\frac 1e
\null^*\Delta_i\null^*\Delta_j f(\null^*v)f(\null^*w)\, dw\,
d\sigma.
  \ee
Different expressions of the operator \fer{corr} can be recovered
owing to the definition of $\Delta$. For the purposes of the
present paper, however, we simply remark that, choosing the test
function $\varphi(v) = |v|^2$, direct computations show that the
correction $D_e(f,f)$ is such that
 \begin{align}\label{temp}
(|v|^2\, , \,D_e(f,f)) \,\, &= \frac {\beta^2}{64\pi}
\int_{\Realr^3}\int_{\Realr^3}\int_{S^2}\left||u|\sigma-u\right|^2
f(v)f(w)dv\, dw\, d\sigma\nonumber \\
&= \frac {1}{4}\beta^2\left[  \int_{\Realr^3}|v|^2f(v)\, dv -
\left|\int_{\Realr^3}v\,f(v)\, dv\right|^2\right].
\end{align}
This reveals  the fundamental fact that the diffusion operator produces a
growth of the second moment proportional to the second moment itself. This
action is clearly different from the action of a linear diffusion operator (a
thermal bath), which induces a growth of the second moment proportional to the
mass. This supports the fact that the Boltzmann equation \fer{BE2} can
produce fat tails.

\subsection{Simple economy market modelling}
\label{eco1}

In one dimension of the "velocity" variable, a similar construction
leads to kinetic models for wealth redistribution \cite{CPT,
MaTo06}. In this case, the variable
 $v \in \Realr_+$ represents
the wealth of the agents, binary collisions are trades between
agents, and the (eventual) power-like tails of the steady
distribution of wealth are known in the pertinent literature as
Pareto tails. Due to the fact that the variable is in $\Realr_+$,
the possible conserved quantities reduce to mass and momentum. In
\cite{CPT} the collision mechanism is given by
\begin{equation}
\label{mixing}
v' = (1-\lambda)v +\lambda w + \eta v; \qquad
 w' = \lambda v + (1-\lambda) w + \eta^* w
\end{equation}
where $0\leq \lambda\leq 1$ represents the constant saving rate
and $\eta$ and $\eta^*$ are random variables with law given by a measure
$\mu(s)$ of zero mean, variance $\beta^2$ and support in
$[-\lambda,+\infty)$. In this way, for all realizations of the
random variable we have $\eta\geq -\lambda$ and wealths after
trading are well defined
i.e.,
remain nonnegative. This is the so-called {\em no debt} condition.
In this
context, the Boltzmann equation \fer{BE2} is replaced by
  \be
(\varphi,\tQ_\lambda(f,f))= \left\langle \int_{\Realr_+}
\int_{\Realr_+} f(v) f(w) \Big[ \varphi (v')- \varphi (v) \Big]
dv\,dw\ \right\rangle. \label{Qweak2}
 \end{equation}
 Here, we use the notation
$$
\left< h\right> :=\int_{-\lambda}^\infty h(s)d\mu(s).
$$
The unique possible collision invariants of the one-dimensional Boltzmann
equation are obtained for $\varphi(v) =1$ and $\varphi(v) =v$.


The weak formulation of the Boltzmann equation can also be rewritten
\begin{equation}
\label{Max1d}
\int_{\rr^+} \varphi(v)\,\tQ_\lambda(f,f)\, dv \!=\!  \frac12
\int_{\rr^+} \!\int_{\rr^+}\!\! f(v) f(w) \left< \varphi (v') +
\varphi (w')- \varphi (v)- \varphi (w) \right> \,dv\, dw .
 \end{equation}
In \eqref{Max1d} the wealth variables $v,w$ are nonnegative quantities, and the
collision mechanism is given by \fer{mixing}.
A one-dimensional Boltzmann type equation of the form
\begin{equation}\label{ibe1d}
 \frac{\p f}{\p t} = \tQ_{\lambda}(f,f)
\end{equation}
based on the binary interaction given in \eqref{mixing} has been
considered in \cite{CPT,MaTo06} and we refer to them for a deeper
discussion of the model. Without loss of generality, we can fix
the initial density $f_0(v) \in {\cal P}_2(\rr)$, with the
normalization condition
\begin{equation}\label{norm1}
  m(t):=\int_{\rr^+} v f(v,t)\, dv =  \bar{m},
\end{equation}
since by choosing $\varphi(v) = v$, \eqref{Max1d} shows that
$m(t)=m(0)$ for all $t\geq0$.\\

As in section \ref{formal}, one splits the collision mechanisms into
a deterministic inelastic part and the random part :
$$v' = v^* + \eta v; \qquad
w' = w^* + \eta^* w
$$
where $v^*, w^*$ are deterministic wealth (corresponding to inelastic
collision with constant restitution coefficient $(1-\lambda)$)
$$
v^* = (1-\lambda)v +\lambda w ; \qquad
w^* = \lambda v + (1-\lambda) w .
$$

A formal Taylor expansion similar to \fer{formal}, in the limit for
$\lambda$ and $\eta$ small, leads to a drift term for the
difference between $(v,w)$ and $(v^*,w^*)$ and a diffusion
term proportional to the variance $\beta^2$.
$$
\varphi(v') = \varphi(v^*) + \eta v  \partial_v \varphi(v^*)  +
\frac 12 \eta^2 v^2  \partial_v^2 \varphi"(v^*)+ \dots
$$
Taking the average
$$
<\varphi(v')> = \varphi(v^*)  + \frac 12 \beta^2 v^2  \partial_v^2 \varphi(v^*)+ \dots,
$$
and on the other hand, the deterministic part gives
$$
\varphi(v^*) =  \varphi(v) + \lambda (w-v)  \partial_v \varphi(v)  +
\frac 12 \lambda^2 (v-w)^2  \partial_v^2 \varphi(v)+ \dots
$$
Inserting these expansions into the weak formulation
of the Boltzmann equation \fer{Qweak2} and rescaling the time
gives
$$
< (\varphi,\tQ_\lambda(f,f))> = \int_{\Realr_+^2}
 f(v) f(w) \Big[ \lambda (w-v) \partial_v \varphi (v)
+ \frac 12 (  \lambda^2 (v-w)^2 + \beta^2 v^2 ) \partial_v^2 \varphi (v) \Big]
dv\, dw.
$$

More precisely, the asymptotics of the
one-dimensional Boltzmann equation for
wealth distribution \fer{Qweak2} for
$\lambda$ sufficiently small,
and in the limit $\frac \lambda{\beta^2} \to \gamma$,
has been studied in \cite{CPT}. In this so-called
"continuous trading limit",
it is proved that the solution
to the Boltzmann equation converges toward the solution to the
Fokker-Planck equation
 \be\label{FP2b}
 \frac{\partial f}{\partial t} = \frac \gamma{2}\frac{\partial^2 }
 {\partial v^2}\left( v^2 f\right) + \frac{\partial }{\partial v}\left((v - \bar m)
 f\right),
 \ee
 which admits a unique stationary state of unit mass, given by the
 $\Gamma$-distribution
 \be\label{equi}
M_\lambda(v)=\frac{(\mu-1)^\mu}{\Gamma(\mu)}\frac{\exp\left(-\frac{\mu-1}{v}\right)}{v^{1+\mu}}
 \ee
  where
  $$ \mu = 1 + \frac{2}{\lambda} >1.
$$
This stationary distribution exhibits a Pareto power law tail for large
velocities. We remark that in \fer{FP2b} the growth of the second moment
follows the same law as the Boltzmann equation \fer{BE2}.


\section{Quick overview of probability metrics}
\label{quick}

In this section, we first briefly recall the main definitions and results
about probability metrics and, more precisely, on Wasserstein ($W_2$)
and Fourier ($d_s$) distances between two probability measures.

\subsection{Wasserstein distances}
\label{Wass1}

Given two probability measures $f,g\in\Pn$, the Euclidean
Wasserstein Distance is defined as
\begin{equation}\label{w2def1}
W_2(f, g) = \inf_{\Pi\in\Gamma} \left\{ \iint_{\Rn \times \Rn}
\vert v -x \vert^2 \, d\Pi(v, x) \right\}^{1/2}
\end{equation}
where $\Pi$ runs over the set of transference plans $\Gamma$, that
is, the set of joint probability measures on $\Rn \times \Rn$ with
marginals $f$ and $g\in \Pn$. From a probabilistic point of view,
the Wasserstein distance can be alternatively defined as
\begin{equation}\label{w2def2}
W_2(f, g) = \inf_{(V,X)\in\tilde\Gamma} \left\{ \expec\left[ \vert
V-X \vert^2 \right] \right\}^{1/2}
\end{equation}
where $\tilde\Gamma$ is the set of all possible couples of random
variables $(V,X)$ with $f$ and $g$ as respective laws. Let us
remark that $W_2$ is finite for any two probability measures with
finite second moments $f,g \in\Ptwo$.

The main properties of the Euclidean Wasserstein distance $W_2$
are summarized in the following proposition. We refer to
\cite{Bolley,Vil03,Vil07} for the proofs and further information
on the connections to optimal mass transport theory.

\begin{proposition}[$W_2$-properties]\label{w2properties}
The space $(\Ptwo,W_2)$ is a complete metric space. Moreover, the
fol\-lowing properties of the distance $W_2$ hold:
\begin{enumerate}
\item[i)] {\bf Optimal transference plan:} The infimum in the
definition of the distance $W_2$ is achieved at a joint
probability measure $\Pi_o$ called an optimal transference plan
satisfying:
$$
W_2^2(f, g) = \iint_{\Rn \times \Rn} \vert v -x \vert^2 \,
d\Pi_o(v, x).
$$

\item[ii)] {\bf Convergence of measures:} Given $\{f_n\}_{n\ge 1}$
and $f$ in $\Ptwo$, the following three assertions are equivalent:
\begin{itemize}
\item[a)] $W_2(f_n, f)$ tends to $0$ as $n$ goes to infinity.

\item[b)] $f_n$ tends to $f$ weakly-* as a measure and
$$
 \int_{\Rn} \vert v \vert^2 \, f_n(v) \, dv \to
\int_{\Rn} \vert v \vert^2 \, f(v) \, dv  \, \mbox{ as } \,
\mbox{n} \to + \infty.
$$
\end{itemize}

\item[iii)] {\bf Convexity:} Given $f_1$, $f_2$, $g_1$ and $g_2$
in $\Ptwo$ and $\alpha$ in $[0,1]$,
$$
W_2^2(\alpha f_1 + (1-\alpha) f_2,\alpha g_1 + (1-\alpha) g_2)
\leq \alpha W_2^2(f_1,g_1) + (1-\alpha) W_2^2(f_2,g_2).
$$
As a simple consequence, given $f,g$ and $h$ in $\Ptwo$,
$$
W_2(h * f,h * g) \leq W_2(f,g)
$$
where $*$ stands for the convolution in $\Rn$.

\item[iv)] {\bf Additivity with respect to convolution:} Given
$f_1$, $f_2$, $g_1$ and $g_2$ in ${\cal P}_2(\Rn)$ with equal mean values,
$$
W_2^2(f_1 * f_2, g_1 * g_2) \leq W_2^2(f_1,g_1) + W_2^2(f_2,g_2).
$$

\end{enumerate}
\end{proposition}

\subsection{Fourier metrics}
\label{Four1}

Given $f\in\Pn$, its Fourier transform or characteristic function
is defined as
\[
\hat{f}(k) = \int_{\Rn} e^{-iv\cdot k}\, df(v).
\]
Given any $s>0$, the Fourier-based metric $d_s$  is defined as
\begin{equation}\label{d2def}
d_s(f,g)= \sup_{k \in \Rno} \frac{|\hat{f}(k)-\hat{g}(k)|}{|k|^s}
\end{equation}
where $\Rno=\Rn-\{0\}$, for any pair of probability measures
$f,g\in\Pn$. This metric was introduced in
\cite{Gabetta-Toscani-Wennberg} and further used in
\cite{CGT,CCG1,Toscani-Villani,GJT}. Only recently, various
applications to the large-time behavior of the dissipative
Boltzmann equation \cite{PT03, BisiCT, BisiCT2} have revealed the
importance of this distance. 
We refer to
\cite{CTPE} for a complete survey of this metric and the proofs of
the statements below.

The metric $d_s$ with $s>0$ is well-defined and finite for any two
probability measures $f,g\in {\cal P}_s(\Rn)$ with
equal moments up to $[s]$ if $s\notin\N$, or equal moments up to
$s-1$ if $s\in \N$. The main properties of the $d_s$ metrics
relevant to the ongoing discussion are summarized in the following
result:

\begin{proposition}\label{dsproperties}
The distances $d_s$ with $s>0$ verify the following properties:
\begin{enumerate}

\item[i)] {\bf Convexity:} Given $f_1$, $f_2$, $g_1$ and $g_2$ in
${\cal P}_s(\Rn)$ with equal moments up to $[s]$ if
$s\notin\N$, or equal moments up to $s-1$ if $s\in \N$ and
$\alpha$ in $[0,1]$,
$$
d_s(\alpha f_1 + (1-\alpha) f_2,\alpha g_1 + (1-\alpha) g_2) \leq
\alpha d_s(f_1,g_1) + (1-\alpha) d_s(f_2,g_2).
$$

\item[ii)] {\bf Superadditivity with respect to convolution:}
Given $f_1$, $f_2$, $g_1$ and $g_2$ in ${\cal P}_s(\Rn)$
with equal moments up to $[s]$ if $s\notin\N$, or equal
moments up to $s-1$ if $s\in \N$,
$$
d_s(f_1 * f_2, g_1 * g_2) \leq d_s(f_1,g_1) + d_s(f_2,g_2).
$$
\end{enumerate}
\end{proposition}


\section{ Large time behavior for economy model}
\label{eco2}

\subsection{Evolution of Wasserstein distance}

The Boltzmann equation \eqref{ibe1d} can be rewritten as
$$
\frac{\p f}{\p t}= \left<f_{p+\eta}*f_q\right> -f ,
$$
where we use the shorthand $f_p(v) = (1/p)f(v/p)$ with
$p=\lambda$ and $q=1-\lambda$. Here, $f$ is extended by 0 to the
whole of $\rr$ in the convolution. The gain operator is defined as
the measure given by
\[
 (\varphi,\tQ_{\lambda}^+(f,f))= \left<\int_{\rr^+} \int_{\rr^+} f(v) \, f(w) \, (\varphi, \delta_{(p+\eta)v+qw}) \, dv \,
 dw\right>
\]
where $\delta_{(p+\eta)v+qw}$ is the Delta Dirac at the
post-collisional velocity $v'$ and $(\cdot,\cdot)$ is the duality
pair between continuous functions and probability measures. In
probabilistic terms, the gain operator is defined as an
expectation:
$$
 \tQ_{\lambda}^+(f,f)= <f_{p+\eta}*f_q> =\expec \left[ \delta_{(p+\eta)V+qW}\right]
$$
where $V$ and $W$ are independent random variables with law $f$
and independent with respect to the random variable $\eta$. Here
the expectation is taken with respect to all random variables.

Let us take two independent pairs of random variables $(V,X)$ and
$(W,Y)$ such that $V$ and $W$ have law $f_1$ and $X$ and $Y$ have
law $f_2$. From the convexity of $W_2^2$ and the independence of
the pairs, it follows that
$$
W_2^2(\tQ_{\lambda}^+(f_1,f_1),\tQ_{\lambda}^+(f_2,f_2))\leq
\expec\left[
W_2^2(\delta_{(p+\eta)V+qW},\delta_{(p+\eta)X+qY})\right]
$$
for any probability densities $f_1,f_2\in{\cal P}_2(\rr)$. Now, the last term
is directly computed as the Euclidean distance between the two points $(p+\eta)V+qW$
and $(p+\eta)X+qY$, and thus,
$$
W_2^2(\tQ_{\lambda}^+(f_1,f_1),\tQ_{\lambda}^+(f_2,f_2))\leq
\expec\left[ |(p+\eta)(V-X)+q(W-Y)|^2 \right].
$$
Using independence of the pairs and taking the pairs to be optimal
couples for the $W_2(f_1,f_2)$ in the probabilistic definition
\eqref{w2def2}, we deduce finally the property
$$
W_2^2(\tQ_{\lambda}^+(f_1,f_1),\tQ_{\lambda}^+(f_2,f_2))\leq
\left[<(p+\eta)^2>+q^2\right] \, W_2^2(f_1,f_2).
$$
Let us define, for $s \ge 1$
 \be\label{key}
 \mathfrak{S}(s) :=  \langle (p +\eta)^s\rangle +q^s  -1;
 \ee
then $\mathfrak{S}(2)=<(p+\eta)^2>+q^2-1= 2\lambda(\lambda
-1)+\beta^2$. It is not difficult to see that the convexity
property of $W_2^2$ together with the Duhamel formula for
\eqref{ibe1d} and the contractive estimate of the gain operator in
$W_2$ leads to the result:

\begin{theorem}\label{cont1d}
Let $f_1(t)$ and $f_2(t)$ be two solutions of the one dimensional Boltzmann
equation \eqref{ibe1d} corresponding to initial values $f_1^0$ and $f_2^0$ in
${\cal  P}_2(\rr^+)$, satisfying conditions \eqref{norm1}. Then, for all times
$t \geq 0$,
\begin{equation}\label{decc}
  W_2(f_1(t), f_2(t)) \leq \exp\left\{ \mathfrak{S}(2) t\right\} W_2(f_1^0,f_2^0).
\end{equation}
If $\beta^2 < 2\lambda(1-\lambda)$, then $\mathfrak{S}(2) <0$, and
the Wasserstein metric decays exponentially to zero in time.
\end{theorem}


\subsection{Evolution of Fourier metrics}

Analogous results for the evolution of the $d_s$-metric
\fer{d2def} have been obtained recently in \cite{MaTo06} by a
suitable generalization of results in \cite{PaTo06}. For the
detailed computations we refer to \cite{MaTo06}. The study of the
evolution of the metric \fer{d2def}, leading to the understanding
of the large-time behavior of the solution to the kinetic equation
\fer{ibe1d}, requires a fine analysis of the quantity \fer{key}. As
shown for the Wasserstein metric in the previous subsection, the
sign of this quantity is in fact related to the contraction
properties of the metric. Moreover, as has been noted in
\cite{MaTo06}, the sign of \fer{key} is also related both to the
number of moments of the solution which remain uniformly bounded
in time, and to the possibility to conclude the existence and
uniqueness of a steady state. The results in \cite{MaTo06} can be
briefly summarized into the following

\begin{theorem}\label{main}
Take $s>0$ with $\mathfrak{S}(s)<\infty$ and let $f_1(t)$ and
$f_2(t)$ be two solutions of the one dimensional Boltzmann
equation \eqref{ibe1d} corresponding to initial values $f_1^0$
and $f_2^0$ in ${\cal P}_r(\rr^+)$, satisfying conditions
\eqref{norm1} with $r = \max\{ 1,s\}$.  Then the following bound
holds:
\begin{equation}\label{conv1}
    d_{s}(f_1(t),f_2(t) ) \leq
    \exp\left\{\mathfrak{S}(s) t\right\}\,d_{s}(f_1^0,f_2^0),
  \end{equation}
  where $ \mathfrak{S}(s)$ is given by \fer{key}.
\end{theorem}

Also, the temporal behavior of the moments is almost completely
determined by the function $\mathfrak{S}(s)$.

\begin{theorem} \label{mom1}
Let $s>1$ and $f_0\in{\cal  P}_s(\rr^+)$ with
$0<\mathfrak{S}(s)<\infty$ and let us denote
$$
M_s^0:=\int_{\rr^+} v^s\,f_0(v)\,dv.
$$
Then, for the weak solution to the Boltzmann equation, the
following estimates hold:
  \begin{enumerate}
  \item If $\mathfrak{S}(s)>0$, then, as $t\to\infty$,
    \begin{align*}
      { \int_{\Realr_+} v^s f(v,t)\,dv } \ge M_s^0 \,{\exp\{\mathfrak{S}(s)t\}
      } + o(1).
    \end{align*}
  \item If $\mathfrak{S}(s)<0$, then the $s$th moment is bounded for all times.
    Moreover, as $t\to\infty$,
    \begin{align*}
       { \int_{\Realr_+} v^s f(v,t)\,dv } \leq M_s^0 \,{\exp\{\mathfrak{S}(s)t\}
      } + o(1).
  \end{align*}
  \end{enumerate}
Here, the remainder terms $o(1)$ converge to zero exponentially
fast.
\end{theorem}

Another important conclusion of the analysis of \cite{MaTo06} is that the
essential function $\mathfrak{S}(s)$ does not only decide whether or not the
steady state $f_\infty$ develops a Pareto tail. In fact, the positive zero of
$\mathfrak{S}(s)$ actually determines the value of the Pareto index.

A comparison of the contraction results for the Boltzmann equation
\eqref{ibe1d} shows that the contraction properties are heavily linked, through
the key function \fer{key}, to the (eventual) formation of tails. While the
situation for equation \eqref{ibe1d} is reasonably well understood, the
corresponding analysis for the Boltzmann equation \fer{BE2} deserves further
investigation. We will discuss equation \fer{BE2} in detail in the
following section.


\section{Large time behavior for stochastic granular media}
\label{seccontr}

Let us consider here the modification of the Inelastic Maxwell Model introduced
in \cite{Bobylev-Carrillo-Gamba}
\begin{align}\label{him}
\frac{\partial f}{\partial t} = \tQ_e(f,f),
\end{align}
where the collision operator is defined weakly as
\begin{align}\label{wfqm2}
(\varphi,\tQ_e(f,f)) & = \frac{1}{4 \pi} \left< \int_{\rr^3} \int_{\rr^3}
\int_{S^2} f(v) f(w) \Big[ \varphi (v')- \varphi (v) \Big] d\sigma \,dv\,dw
\right> .
\end{align}
As discussed in the introduction, the collision mechanism relies
on a random coefficient of restitution,
\begin{eqnarray}
v' &=& \frac12 (v+w) + \frac{1-\et}{4}u+\frac{1+\et}{4} |u|\sigma
\nonumber
\\[-2mm]
\label{colmech2}
\\[-2mm]
w' &=& \frac12 (v+w) - \frac{1-\et}{4}u-\frac{1+\et}{4} |u|\sigma. \nonumber
\end{eqnarray}
As before, we write $u=v-w$, $\et=e+\eta$ and $\eta$ is a
real-valued random variable, with zero mean and variance $\beta^2$,
given by a measure $\mu(s)$ with support on $[-e,\infty)$. Here,
$<\cdot>$ means the expectation with respect to $\eta$, i.e., the
integral over $\rr$ with respect to $\mu$.

It is quite straightforward to check that conservation of mass and
momentum remains and that
$$
\left<|v'|^2 + |w'|^2- |v|^2 - |w|^2 \right> =0
$$
for the model that is conservative in the mean in which $\beta^2=1-e^2$. From
\eqref{wfqm2}, we deduce that the temperature evolution is
$$\frac{d}{d t} \int_{\rr^3} |v|^2 \,f(t,v)\,dv = 0,$$ and thus we
deduce that $\theta(t)=\theta(0)$ for all times $t\geq0$ and we
will fix it to one for convenience.


\subsection{Evolution of Wasserstein distance }

Given a probability measure $f$ on $\rr^3$, the gain operator is
in fact a probability measure $\tQ_e^+(f, f)$ defined by
\[
 (\varphi,\tQ_e^+(f,f))= \left<\int_{\rr^3} \int_{\rr^3} f(v) \, f(w) \, (\varphi, {\cal U}_{v, w,\eta}) \, dv \,
 dw\right>
\]
where ${\cal U}_{v, w,\eta}$ is the uniform probability
distribution on the sphere $S_{v, w}$ with center $c_{v, w} =
\frac{1}{2} (v+w)+ \left[\frac{1-\et}{4}+\eta\right](v-w)$ and
radius $r_{v, w} = \frac{1+\et}{4} |v-w|$ as in
\cite{Bolley-Carrillo}. In probabilistic terms, the gain operator
is defined as an expectation:
$$
\tQ_e^+(f, f)=\expec \left[ {\cal U}_{V, W,\eta}\right]
$$
where $V$ and $W$ are independent random variables with law $f$
and independent of the law of $\eta$. As in
\cite{Bolley-Carrillo}, we get the following result:

\begin{theorem}\label{contrqw2} Given $f$ and $g$ in
$\Ptwot$ with equal mean velocity, then
\[
W_2(\tQ_e^+(f, f) , \tQ_e^+(g, g) ) \leq  W_2(f, g).
\]
\end{theorem}

\proof Let us take two independent pairs of random variables
$(V,X)$ and $(W,Y)$ such that $V$ and $W$ have law $f$ and $X$ and
$Y$ have law $g$. Also, let us take two independent random
variables $\eta$ and $\tilde{\eta}$ with law $\mu$. Convexity of
$W_2^2$ implies
\begin{equation}
W_2^2(\tQ_e^+(f,f), \tQ_e^+(g,g)) = W_2^2(\expec \left[ {\cal
U}_{V, W,\eta}\right], \expec \left[ {\cal U}_{X,
Y,\tilde{\eta}}\right])\leq \expec \left[ W_2^2({\cal U}_{V,
W,\eta}, {\cal U}_{X, Y,\tilde{\eta}})\right]\label{tech8}
\end{equation}
where the expectation is taken with respect to the joint
probability density in $\rr^{14}$ of the six random variables.
Here, the independence of the pairs of random variables has been
used.

As proved in \cite{Bolley-Carrillo}, the $W_2^2$ distance between
the uniform distributions on the sphere with center $O$ and radius
$r$, ${\mathcal U}_{O,r}$, and on the sphere with center $O'$ and
radius $r'$, ${\mathcal U}_{O',r'}$, in $\rr^3$ is bounded by
$\vert O'-O \vert^2 + (r' - r)^2$.

We now estimate the right-hand side of \eqref{tech8} by using the
formulas for the center and radii of the spheres given in
\eqref{colmech2} to deduce
\begin{align*}
W_2^2(\tQ_e^+(f,f), \tQ_e^+(g,g)) \!\leq  \!& \; \left<\frac{5 - 2
\, \et + \et^2}{8}\right> \expec \left[ \vert V-X \vert^2 \right]
+ \, \left<\frac{(1 + \et)^2}{8}\right> \, \expec \left[
\vert W-Y \vert^2 \right] \\
& + \, \left<\frac{1-\et^2}{4}\right> \, \expec \left[ (V-X) \cdot
(W-Y)\right]
\end{align*}
where the Cauchy-Schwartz inequality has been used.

Finally, we take both pairs $(V,X)$ and $(W,Y)$ as independent
pairs of variables with each of them being an optimal couple for
$W_2(f,g)$ in the probabilistic definition \eqref{w2def2}
to obtain
\begin{align*}
W_2^2(\tilde{Q}_e^+(f,f), \tilde{Q}_e^+(g,g)) \leq  & \; \frac{3 +
e^2+\beta^2}{4} \, W_2^2(f,g) + \, \frac{1-e^2-\beta^2}{4} \,
\expec \left[ (V-X) \cdot (W-Y)\right],
\end{align*}
where the last term is zero because the random variables are independent
and have equal means. Since $\beta^2 =1-e^2$
in the conservative case, the result is proved.\endproof

\

As a consequence of the previous property of the gain operator, we
draw the following conclusion about controlling the distance
between any two solutions of \eqref{him} in the conservative case.

\begin{theorem}\label{contrw2eb} If $f_1$ and $f_2$ are two solutions
to \eqref{him} with respective initial data $f_1^0$ and $f_2^0$ in
$\Ptwot$ with zero mean velocity, then, for all $t \geq 0$,
\begin{align*}
W_2^2(f_1(t), f_2(t)) \leq W_2^2(f_1^0, f_2^0).
\end{align*}
\end{theorem}

\proof Duhamel's formula for \eqref{him} reads as
\[
f_i(t) = {\rm e}^{-t} \, f_i^0 + \int_0^{t} {\rm e}^{-(t-s)} \,
\tilde{Q}_e^+(f_i(s), f_i(s)) \, ds, \qquad i=1,2.
\]
As before, the convexity of the squared Wasserstein distance in
Proposition \ref{w2properties} and the contraction of the gain
operator in Theorem \ref{contrqw2} imply
\begin{align*}
W_2^2(f_1(t), f_2(t)) & \leq {\rm e}^{-t} \, W_2^2(f_1^0, f_2^0) +
\int_0^{t} \!\! {\rm e}^{-(t-s)} \, W_2^2 \big(
\tilde{Q}_e^+(f_1(s), f_1(s)), \tilde{Q}_e^+(f_2(s), f_2(s)) \big) \, ds \\
& \leq {\rm e}^{-t} \, W_2^2(f_1^0, f_2^0) + \, \int_0^{t} \!\!
{\rm e}^{-(t-s)} \, W_2^2 (f_1(s), f_2(s))\, ds.
\end{align*}
Therefore, the function $y(t) = {\rm e}^{t} \, W_2^2(f_1(t),
f_2(t))$ satisfies the inequality
 \[
 y(t) \leq y(0) + \int_0^{t} y(s) \, ds
 \]
and thus $y(t) \leq y(0) \, {\rm e}^{t}$ by Gronwall's lemma,
concluding the argument.
\endproof


\subsection{Evolution of Fourier metrics}

We start by writing a closed form of the Boltzmann equation in
Fourier variables. In fact, it is not difficult using Bobylev's
identity in \cite{Boby75,Bobylevid,Bobylev,Bobylev-Carrillo-Gamba}
to get
$$
\widehat{\tQ_e^+(f,f)} = \frac{1}{4 \pi} \left<\int_{S^2} \hat{f}
(t, k_-) \hat{f} (t, k_+) \, d\sigma\right>
$$
where
$$
k_-=\frac{1+\et}{4} k - \frac{1+\et}{4}|k|\sigma \quad \mbox{and}
\quad k_+=\frac{3-\et}{4}\, k + \frac{1+\et}{4}\, |k|\sigma\,.
$$
Let us start by analyzing the evolution of the distance $d_2$ that
in view of the properties in Propositions \ref{w2properties} and
\ref{dsproperties} should verify the same non-strict contraction
as the transport distance $W_2$.

\begin{theorem}\label{contrqd2} Given $f$ and $g$ in $\P_2(\rr^3)$
with equal mean velocity,
\[
d_2(\tQ_e^+(f, f),\tQ_e^+(g, g) ) \leq \frac{3+ e^2+\beta^2}{4} \,
d_2(f,g).
\]
\end{theorem}

\proof Using the Fourier representation formula above, we deduce
$$
\frac{\widehat{\tQ_e^+(f, f)}(k) -
\widehat{\tQ_e^+(g,g)(k)}}{|k|^2} = \frac{1}{4 \pi}\!\left<
\int_{S^2} \left[ \frac{\hat{f}(k_-) \hat{f}(k_+) - \hat{g}(k_-)
\hat{g}(k_+)}{|k|^2} \right] d\sigma\right>
$$
for all $k\in\Rto$. We now estimate the integrand as
\begin{align*}
\left| \frac{\hat{f}(k_-) \hat{f}(k_+)-\hat{g}(k_-)
\hat{g}(k_+)}{|k|^2} \right| &\leq \sup_{k \in \Rto} \left\{
  \frac{|\hat{f}(k)-\hat{g}(k)|}{|k|^2} \right\} \left( \frac{|k_-|^2+
    |k_+|^2}{|k|^2} \right) \\
    & = d_2(f,g) \left( \frac{|k_-|^2+
    |k_+|^2}{|k|^2} \right),
\end{align*}
and thus
$$
d_2(\tQ_e^+(f, f) , \tQ_e^+(g, g) ) \leq \frac{1}{4 \pi}\!
\left<\int_{S^2} \left( \frac{|k_-|^2+
    |k_+|^2}{|k|^2} \right) d\sigma \right> \, d_2(f,g).
$$

We observe that
 \[
 \frac{|k_-|^2+|k_+|^2}{|k|^2}
 \]
  is a function of
the angle between the unit vectors $k/|k|$ and $\sigma$ and the
random variable $\eta$, and that
$$
I:=\frac{1}{4 \pi} \left<\int_{S^2} \frac{|k_-|^2+|k_+|^2}{|k|^2}
\,d\sigma \right>= \frac{3+e^2+\beta^2}{4}.
$$
In fact, we can compute
\begin{equation} \label{modulusk+-}
\begin{split}
&\displaystyle |k_-|^2= |k|^2 \left( \frac{1+\et}{4} \right)^2 2\,
\Big( 1- \cos \vartheta \Big)
\\*[.3cm] &\displaystyle |k_+|^2=
|k|^2 \left[ \left( \frac{3-\et}{4} \right)^2 + \left(
\frac{1+\et}{4} \right)^2 + 2 \left( \frac{3-\et}{4} \right)
\left( \frac{1+\et}{4} \right) \cos \vartheta \right]
\end{split}
\end{equation}
where $\vartheta$ is the angle between the unit vectors $k/|k|$
and $\sigma$ from which the value of $I$ is obtained. Putting
together previous estimates we get the contraction in $d_2$ with
the same constant as $W_2^2$ as desired.
\endproof

\

Now, let us see that we can also control Fourier-based distances
with exponent $2+\alpha$, $\alpha\in [0,\infty)$. Let us set
\begin{align}\label{key2}
\A(\alpha,e,\eta) & := \displaystyle \frac{1}{2}\left< \int_0^\pi
\left\{ \left[ \left( \frac{1+\et}{4} \right)^2 2 (1- \cos
\vartheta ) \right]^{\frac{2+\alpha}{2}}
\right.\right. \vspace{0.2 cm} \nonumber\\
& \,\,\,\,\,\,\left.\left. +\, \left[ \left( \frac{3-\et}{4} \right)^2 + \left(
\frac{1+\et}{4} \right)^2 + 2 \left( \frac{3-\et}{4} \right) \left(
\frac{1+\et}{4} \right) \cos \vartheta \right]^{\frac{2+\alpha}{2}} \right\}
\sin \vartheta\, d\vartheta \right>
 \vspace{0.2 cm} \nonumber\\
& =  \displaystyle \frac{2}{4+\alpha} \left< \left( \frac{1+\et}{2}
\right)^{2+\alpha} + \frac{1- \left| \frac{1-\et}{2}
\right|^{4+\alpha}}{1-\left| \frac{1-\et}{2} \right|^2} \right> .
\end{align}
Whenever there is no confusion, i.e. for $e$ and $\eta$ fixed, we
will denote just by $\A(\alpha)$ the above constant.

\begin{theorem}\label{contrqd2alpha} Given $f,g\in
\P_{2+\alpha}(\rr^3)$ with equal moments up to order $2+[\alpha]$, there
exists an explicit constant $\mathfrak{A}(\alpha,e,\eta)>0$ given by \fer{key2}
such that
\[
d_{2+\alpha}(\tQ_e^+(f, f) , \tQ_e^+(g, g) ) \leq \A(\alpha,e,\eta) \,
d_{2+\alpha}(f,g).
\]
\end{theorem}

\proof As in the proof of the previous theorem, we compute
\begin{align*}
\left| \frac{\widehat{\tQ_e^+(f, f)}(k) -
\widehat{\tQ_e^+(g,g)}(k)}{|k|^{2+\alpha}} \right| & = \frac{1}{4
\pi} \left|\left< \int_{S^2} \frac{\hat{f}(k^+) \hat{f}(k^-)
-\hat{g}(k^+) \hat{g}(k^-)}{|k|^{2+\alpha}}\, d\sigma \right>\right| \vspace{0.1 cm}\nonumber\\
 & \leq A \, \sup_{k \in \Rto} \frac{|\hat{f}(k)-\hat{g}(k)|}{|k|^{2+\alpha}}
\end{align*}
where $A$ is given by
\begin{equation}
A := \frac{1}{4 \pi} \left< \int_{S^2} \frac{|k_+|^{2+\alpha} +
|k_-|^{2+\alpha}}{|k|^{2+\alpha}}\, d\sigma \right>.\label{Aalpha}
\end{equation}
By inserting the expressions of $k_-$ and $k_+$ into
\eqref{Aalpha} and computing the integral we conclude
$A=\A(\alpha,e,\eta)$ and the proof follows.
\endproof

\

As a consequence, we obtain an estimate on contraction/expansion
of the Fourier distances $d_{2+\alpha}$ between solutions.

\begin{theorem}\label{contrdseb} Let $\alpha>0$ be such that
$\mathfrak{A}(\alpha,e,\eta) <\infty$. Let $f_1$ and $f_2$ be two
solutions to~\eqref{him} corresponding to initial values~$f_1^0$,
$f_2^0$ with equal moments up to $2+[\alpha]$. Then, for all $t
\geq 0$,
\begin{equation} \label{d2alphadecayeb}
d_{2+\alpha}(f_1(t),f_2(t)) \leq d_{2+\alpha}(f_1^0,f_2^0)\, {\rm
e}^{-C(\alpha,e,\eta) t},
\end{equation}
with $C(\alpha,e,\eta)=1-\A(\alpha,e,\eta)$.
\end{theorem}

\proof The Fourier expression of equation \eqref{him} is given by
$$
\frac{\p \hat{f}}{\p t} = \frac{1}{4 \pi} \int_{S^2} \hat{f}(k_+)
\hat{f}(k_-) d\sigma - \hat{f} = \widehat{\tilde{Q}_e^+(f,f)} -
\hat{f},
$$
whose solution satisfies
\begin{equation}
\hat{f}(t,k) = {\rm e}^{-t} \hat{f}(0,k) + \int_0^t {\rm e}^{-
(t-s)} \widehat{\tilde{Q}_e^+(f,f)}(s,k) \,ds\,. \label{finalg}
\end{equation}

Taking the expressions of the two solutions $\hat{f}_1(t)$ and
$\hat{f}_2(t)$ in \eqref{finalg}, subtracting them and dividing by
$|k|^{2+\alpha}$ with $k\in\Rto$, we get
$$
{\rm e}^{t} \frac{(\hat{f}_1-\hat{f}_2) (t, k)}{|k|^{2+\alpha}} =
\, \frac{\hat{f}_1(0,k)-\hat{f}_2(0,k)}{|k|^{2+\alpha}}\, +
\int_0^t \!{\rm e}^{s} \frac{\Big(
\widehat{\tilde{Q}_e^+(f_1,f_1)}- \widehat{\tilde{Q}_e^+(f_2,f_2)}
\Big) (s,k)}{|k|^{2+\alpha}}\, ds.
$$
Using Theorem \ref{contrqd2alpha} and taking the supremum in
$k\in\Rto$, we obtain
\begin{eqnarray*}
{\rm e}^{t} d_{2+\alpha}(\hat{f}_1,\hat{f}_2)(t) \leq d_{2+\alpha}
\big(\hat{f}_1(0),\hat{f}_2(0) \big) + \A(\alpha,e,\eta) \int_0^t
\!\!{\rm e}^{s} d_{2+\alpha}(\hat{f}_1,\hat{f}_2)(s)ds.
\end{eqnarray*}
Let us set $w(\tau)={\rm e}^{t} d_{2+\alpha}
(\hat{f}_1,\hat{f}_2)(t)$. Then
$$
w(t) \leq w(0) + \A(\alpha,e,\eta)\int_0^t w(s)\, ds ,
$$
which by Gronwall's inequality implies $w(t) \leq w(0)\, {\rm
e}^{\A(\alpha,e,\eta) t}$, concluding the proof.
\endproof

\

The function $\A(\alpha):[0,\infty)\longrightarrow\rr^+$ is convex
by direct inspection. Taking into account that $\mathfrak{A}(0) =
1$,  there are only three possible scenarios for the qualitative
behavior of $\mathfrak{A}$. These are characterized by the sign of
$\mathfrak{A}'(0)$. In case $\mathfrak{A}'(0) \geq 0$, the
function $\mathfrak{A}(\alpha)$ has a minimum at $\alpha=0$ due
to convexity, and thus $\mathfrak{A}(\alpha) > 1$ for all
$\alpha>0$. In this case, there does not exist any $\bar \alpha \in
\Realr_+$ such that $\mathfrak{A}(\bar \alpha)<1$ and there are no
contraction, only expansion, estimates of $d_s$ for $s> 2$.

Suppose that $\mathfrak{A}'(0) <  0$. In this case, the
contraction properties of $d_s$ depend on whether
 \[
 \lim_{\alpha \to \infty} \mathfrak{A}(\alpha) < 1
 \]
 or
 \[
 \lim_{\alpha \to \infty} \mathfrak{A}(\alpha)  > 1.
 \]
In the former case, $\mathfrak{A}(\alpha)< 1$ for $\alpha >0 $.
Theorem \ref{contrqd2alpha} then implies that the $d_s$-metric is
contractive for all values of the parameter $s>2$. In the latter,
since $\mathfrak{A}(0) = 0$,
the convex function
$\mathfrak{A}(\alpha)$ has a minimum attained at some point
$\tilde \alpha
> 0$, and at the same time there exists $\bar \alpha > \tilde \alpha$ for which
$\mathfrak{A}(\bar \alpha)= 1$. Thus, $\mathfrak{A}(\alpha)< 1$ in
the interval $0 < \alpha < \bar \alpha $, and at the same time
$\mathfrak{A}(\alpha)> 1$ for $\alpha >\bar \alpha$. In this case
Theorem \ref{contrqd2alpha} implies that the Boltzmann equation is
contractive up to but not including order $\bar \alpha$.

\begin{remark}
In order to clarify the behavior of $\A(\alpha,e,\eta)$, we can
fix the random variable $\eta$ to assume only two values, while
respecting conditions \fer{re}. This can be done by assuming that
$\eta$ only takes the value ${\sqrt{1-e^2}}/\varrho$ with
probability ${\varrho^2}/(1+\varrho^2)$ and the value
${\sqrt{1-e^2}}\varrho$ with probability ${1}/(1+\varrho^2)$. By
varying the parameters $\varrho$ and $e$ one encounters the whole
variety of possible behaviors of the function $\A(\alpha,e,\eta)$.
Since
 \[
\A(\alpha,e,\eta) =  \displaystyle \frac{2}{4+\alpha} \left< \left(
\frac{1+\et}{2} \right)^{2+\alpha} + \frac{1- \left| \frac{1-\et}{2}
\right|^{4+\alpha}}{1-\left| \frac{1-\et}{2} \right|^2} \right>,
 \]
$\A(\alpha,e,\eta)$ results in the sum of four contributions, one
of which is
 \[
C(\alpha,e,\eta) = \displaystyle
\frac{1}{1+\varrho^2}\frac{2}{4+\alpha}\left( \frac{1+e +
\sqrt{1-e^2}{\varrho}}{2} \right)^{2+\alpha}.
 \]
For any fixed values of $\bar\alpha >0$ and $e$, since the
numerator grows like $\varrho^{2+\alpha}$,  we can choose
$\varrho>>1$ in such a way that $C(\alpha,e,\eta) >1$, and Theorem
{\rm\ref{contrqd2alpha}} implies that the Boltzmann equation is
contractive up to but not including order $\bar \alpha$.

On the other hand, choosing for example $\alpha= 2$ to simplify computations, one
obtains easily
 \be\label{A4}
\A(2,e,\eta) = \displaystyle \frac{1}{3} \left< \left( \frac{1+\et}{2}
\right)^{4} + 1 + \left( \frac{1-\et}{2} \right)^{2} + \left( \frac{1-\et}{2}
\right)^{4}\right> =
  \frac{23 - e + \langle \et^4 \rangle}{24}.
 \ee
Choosing now $1-e <<1$, and $\varrho = \sqrt{1- e^2}/e$, one
obtains that $\et$ assumes the value $0$ with probability $1-
e^2$ and the value $1/e$ with probability $e^2$. Therefore
$\langle \et^4 \rangle = 1/e^2$, which implies $\A(2,e,\eta) < 1$
as long as $1/e^2 -e < 2$. In this second case Theorem
{\rm\ref{contrqd2alpha}} implies that the Boltzmann equation is
contractive at least up to order $4$.
\end{remark}

\subsection{Existence and uniqueness of regular isotropic steady states}

Existence and uniqueness of steady states, as well as the size of their
overpopulated tails, can be derived in full generality (that is, without imposing
restrictive conditions on the random coefficient of restitution) by adapting
to the present situation the methodology of \cite{BisiCT2}, which refers to the
inelastic Boltzmann equation for Maxwell molecules. This methodology, in fact,
is based only on the contractivity properties of the $d_s$-metric, which are
analogous to Theorems \ref{contrqd2alpha} and \ref{contrdseb}.

It has to be remarked that the approach in \cite{BisiCT2} is not suitable to
recover the (eventual) regularity of the steady profile. A regularity result
for the steady state of the inelastic Boltzmann equation for Maxwell molecules
has been obtained in a recent paper by Bobylev and Cercignani
\cite{Bobylev-Cercignani2}. In this paper they were concerned with properties
of the self-similar profiles of the Boltzmann equation for both elastic and
inelastic collisions, and, in addition to the existence, they obtained  results
on the regularity of the steady profiles by showing that the Fourier transform
of the steady profile satisfies a suitable upper bound. Their method takes
advantage of the existence of a super-solution  to the rescaled equation in
Fourier variables (BKW-mode). In our collisional setting, the situation is more
involved, and it requires a precise analysis.

In Fourier variables, the steady state of \fer{BE2} is a solution of the
integral equation
 \be\label{sta1}
\frac{1}{4 \pi} \left<\int_{S^2} \hat{f} (k_-) \hat{f} (k_+) \, d\sigma\right>
= \hat{f}(k),
 \ee
where $k_+$ and $k_-$ are given by the relations
$$
k_-=\frac{1+\et}{4} k - \frac{1+\et}{4}|k|\sigma \quad \mbox{and} \quad
k_+=\frac{3-\et}{4}\, k + \frac{1+\et}{4}\, |k|\sigma\,.
$$
Since isotropy is not destroyed by the collision operator, by choosing
isotropic initial values, one concludes with the isotropy of the (eventual)
steady state. Taking this property into account, the following result can be
obtained as a consequence of Theorem \ref{contrdseb} (see \cite{BisiCT2} for
details).

\begin{corollary}\label{exissteadyebds}
Equation~\eqref{him} has a unique isotropic  steady state $f_\infty$ in the set
of isotropic probability measures with unit mass, zero mean velocity and unit
temperature. Moreover, given any solution $f$ to \eqref{him} for the initial
data $f_0\in \P_2(\rr^3)$ with zero mean velocity and unit pressure tensor,
$$
d_{2+\alpha}(f(t),f_\infty) \leq d_{2+\alpha}(f_0,f_\infty)\, {\rm
e}^{-C(\alpha,e,\eta) t}
$$
for all $t \geq 0$, $0<\alpha<1$. Thus, if $\A(\alpha,e,\eta) < 1$, $f(t)$
converges to the stationary state as $t\to\infty$ in the $d_{2+\alpha}$ sense.
\end{corollary}

\begin{remark}
The previous result shows that the stationary states attract all solutions with
initial data having zero mean velocity and unit pressure tensor. The assumption
of having unit pressure tensor can be weakened to having initial unit temperature
by proceeding similarly to the homogeneous cooling state analysis in
{\rm\cite{Bobylev-Cercignani-Toscani,BisiCT2}}.
\end{remark}

Let us define
\begin{equation} \label{+-}
\begin{split}
&\displaystyle a^2(e,\eta,\theta)= \frac{|k_-|^2}{|k|^2}= \left(
\frac{1+\et}{4} \right)^2 2\, \Big( 1- \cos \vartheta \Big)
\\*[.3cm] &\displaystyle b^2(e,\eta,\theta)= \frac{|k_+|^2}{|k|^2}=
 \left[ \left( \frac{3-\et}{4} \right)^2 + \left( \frac{1+\et}{4}
\right)^2 + 2 \left( \frac{3-\et}{4} \right) \left( \frac{1+\et}{4} \right)
\cos \vartheta \right]
\end{split}
\end{equation}
 Recalling the definition of $k_+$ and $k_-$ given in
\fer{modulusk+-}, it is immediate to show that
 \be\label{cc1}
 a+b \ge 1; \qquad  \frac{1}{2} \left<\int_{0}^\pi (a^2 +b^2) \,
 \sin \theta d\theta\right> =1
 \ee
The first property in \fer{cc1} is a direct consequence of the equality $k_+
+k_- = k$, while the second is the equality $\mathfrak{A}(0) =1$ in
\fer{Aalpha}.
 Let us set $x= |k|$. Then, for any function $\psi(x)$, the
 left-hand side of \fer{sta1} can be rewritten
 in the form
 \be\label{sta2}
 R[\psi(x)] = \frac{1}{2} \left<\int_{0}^\pi \psi(ax) \psi(bx) \,
 \sin \theta d\theta\right>.
 \ee
Under the conditions of Corollary \ref{exissteadyebds}, the Boltzmann equation
has a unique steady state $\ff_\infty(x)$, of unit mass, zero mean velocity and
unit second moment.

Let us remark that $0 \le R[\psi] \le 1$ if $0 \le \psi \le 1$, and $R[\psi]
\le R[\phi]$ if $0 \le \psi \le \phi$. Hence the iteration is monotone
increasing and converges point-wise if we choose the initial approximation $0
\le  \phi_0 \le 1$ in such a way that $\phi_0 \le R[\phi_0]$.
 As observed in  \cite{Bobylev-Cercignani2}, $\phi_0(x) = \exp\{-x^2/2\}$
allows us to obtain a monotone increasing sequence.  In fact, since the
 function $e^{-r}$, $r \ge 0$ is convex, by Jensen's inequality we obtain
 \be\label{in5}
 \left\langle {\rm e}^{-\frac 12(a^2+b^2)x^2}\right\rangle \ge e^{ -\left\langle
\frac 12(a^2+b^2)x^2 \right\rangle } = e^{-{x^2}/2}.
 \ee
This implies that the limit $\ff_\infty(x) \ge 0$.  The trivial limit
$\ff_\infty(x)=1$ will be excluded if there exists a non-zero function
$\phi_0(x)$ such that
 \be\label{in4}
 \phi_0(x) \le \psi_0(x),
 \ee
and at the same time $\psi_0(x)$ generates a monotone decreasing sequence.

Inspired by the ideas of Desvillettes et al in \cite{DFT}, given a fixed
positive constant $\R$, we introduce the fixed point operator
\begin{align*}
  \Rop[\psi](x) &:= \left\{ \begin{array}{cl}
      \ff_\infty(x) & \mbox{if $x< \R$} \\
      R[\psi(x)] & \mbox{if $x \geq \R$}
      \end{array} \right.
\end{align*}
on bounded complex functions $\psi:\setR\to\setC$. Notice that $\Rop$ is
closely related to the Fourier transform of the collision kernel.

\begin{lemma}\label{bound}
Let $0 \le \ff_\infty(x)\le 1$ be the steady state of the Boltzmann equation,
and let us define
\begin{align*}
  \psi_0(x) &:= \left\{ \begin{array}{cl}
      \ff_\infty(x) & \mbox{if $ x < \R$} \\
       \exp(-\mu x)& \mbox{if $ x \geq \R$.}
      \end{array} \right.
\end{align*}
Then, if the random variables $a(e, \eta, \theta)$ and $b(e, \eta, \theta)$ are
such that
 \be \label{condd}
 P(a < \delta) + P(b < \delta) \to 0 \qquad {\rm as} \,\,\, \delta \to 0 ,
 \ee
 there exist positive constants  $\R$ and $\mu$ such that
  \[
\Rop[\psi_0](x) \le \psi_0(x).
 \]
 \end{lemma}

 \proof
Since the steady state is an isotropic probability density function  of unit
mass, zero mean velocity and unit second moment, there exist positive constants
$M < 1/2$ and $\R$ such that (cfr. \cite{DFT})
 \be \label{max}
 0 \le \ff_\infty(x) \le e^{-Mx^2} \qquad {\rm if}\,\,\, x \le \R.
 \ee
Hence, we can fix $\R$ and $M$ to obtain
 \be\label{bbb}
 \psi_0(x) \le e^{-Mx^2} \qquad {\rm if}\,\,\, x \le \R.
 \ee
Clearly, thanks to the definition of $\psi_0$, if $x \le \R$, there is nothing
to prove. Therefore, let us consider the possible cases corresponding to $x
>\R$. Since $a+b \ge 1$, if both $ax \ge \R$, $bx \ge \R$,
\[
 \langle \psi_0(ax)\psi_0(bx)e^{\mu x} \rangle \le 1.
 \]
If now both $ax <\R$ and $bx <\R$, using bound \fer{bbb}, we obtain
 \[
 \langle \psi_0(ax)\psi_0(bx)e^{\mu x} \rangle \le \langle e^{g(x)} \rangle,
 \]
 where
 \[
 g(x) = \mu x - M(a^2 + b^2)x^2.
 \]
 Since $ a +b \ge 1$, it follows that $a^2 + b^2 \ge 1/2$. Thus
  \be \label{b2}
  g(x) \le \mu x - \frac 12 (a^2 + b^2)x^2 \le 0 \qquad {\rm if} \,\,\, \mu \le \frac
  M{2\R}.
  \ee
 Consider now the case in which $ax \le  \R$, while $bx >\R$. In this case
\[
 \langle \psi_0(ax)\psi_0(bx)e^{\mu x} \rangle \le \langle e^{h(x)} \rangle,
 \]
 where
 \[
 h(x) = \mu(1-a) x - M b^2 x^2.
 \]
Since $ a +b \ge 1$, it follows that $b \ge 1-a$, and
 \be\label{b3}
  h(x) \le z(bx) = \mu bx - M(bx)^2 \le \frac{\mu^2}{4M^2}.
   \ee
In fact, the function $z(r)$ has a maximum at $\bar r = \mu/(2M)$. Moreover,
since $z(r)$ decreases for $r >\bar r$, if $ r \ge 3 \bar r$,
 \be\label{b4}
 z(r) \le z(3\bar r) = -3 \frac{\mu^2}{4M^2} .
 \ee
Let us split the calculation of the mean value into the sets $A = \{ bx \ge 3
\bar r\}$ and $A^c=  \{ bx < 3 \bar r\}$. Thanks to  conditions \fer{b3} and
\fer{b4} one obtains
 \be\label{mm}
\langle e^{h(x)}\rangle \le P(A) \exp\left\{-3 \frac{\mu^2}{4M^2}\right\} +
P(A^c)\exp\left\{ \frac{\mu^2}{4M^2}\right\}.
 \ee
Let us set $\delta = 3\bar r = 3\mu/(2M)$. By hypothesis, since $x >\rho$,
 \[
 P(A^c) = P(bx < \delta) \le P(b\rho < \delta)  \to 0 \qquad {\rm if} \,\,\, \delta \to 0.
 \]
 Consider that we can rewrite \fer{mm} as
 \[
\langle e^{h(x)}\rangle \le (1- P(A^c) \exp\left\{- \frac 13 M\delta^2 \right\}
+ P(A^c)\exp\left\{  \frac 13 M\delta^2\right\} =
 \]
 \[
 1 - \frac 13 (1- 4P(A^c)) M \delta^2 + o(\delta^2) \le 1
 \]
 if $\delta$ is sufficiently small. Now, this condition on $\delta$ can be
 satisfied by choosing $\mu$ sufficiently small. This is not in contrast with
condition \fer{b2}. Since the case in which $ax >  \R$ while $bx \le \R$ can
be treated likewise, the lemma is proven.
\endproof
 \

\begin{remark}
Condition \fer{condd} excludes some pathological situations related to the
definition of the random variable $\eta$ \which describes the randomness of the
coefficient of restitution $e$. For example, condition \fer{condd} is violated
if $\eta$ is concentrated on some particular point,
 \[
 P( \eta = 1-e) = p >0 .
 \]
In this case, in fact, $P(b(e,\eta, 0) = 0) = p$, and condition \fer{condd} is
false.
\end{remark}

Lemma \ref{bound} implies that, starting from $\psi_0$, the iteration process
leads to a monotone decreasing sequence. On the other hand, it is clear that,
for $\mu$ sufficiently small,
 \[
 0 \le \phi_0(x) \le \psi_0(x) \le 1.
 \]

Given $\mu>0$, define $K_\mu$ as the set of functions $\psi$ with $\psi(0)=1$,
$\psi'(0)=\ff_\infty'(0)$, and satisfying the estimates
\begin{align}
  \label{eq.gevrey}
  |\psi(x)| \leq \exp(-\kappa x^2) \quad \mbox{for $ x < \R$},
  \qquad
  |\psi(x)| \leq \exp(-\mu x) \quad \mbox{for $x \geq \R$}.
\end{align}

The previous inequalities prove the following

\begin{theorem}
For any pair of functions $a$ and $b$ satisfying conditions \fer{cc1}and
\fer{condd}, the integral equation \fer{sta2} has a nontrivial solution
$\ff_\infty(x)$ such that $\ff_\infty(x)$ belongs to the Gevrey class $K_\mu$
defined by \fer{eq.gevrey}.
\end{theorem}

\begin{remark}
An analogous regularity result can be proven for the steady state to the
one-dimensional kinetic model \fer{ibe1d} \cite{MaTo06}. In this case, it is
important to know that the mean wealth of the stationary state is equal to one.
\end{remark}

\subsection{Fat tails of stationary states}

In this work, we will only examine the case of the fourth moment,
postponing the complete analysis of moment evolution to future
research. Here, we will show that under certain conditions on the
random variable, the fourth moment diverges or is controlled
uniformly.

\begin{lemma}\label{prop-mom4} Let the restitution coefficient $e$ and
the random variable $\eta$ be chosen so that $\A(2,e,\eta) < 1$. If
$f^0$ is a Borel probability measure on $\rr^3$ such that
$$
\int_{\rr^3} \vert v \vert^4 \, f^0(v) \, dv<\infty,
$$
then the solution $f$ to \eqref{him} with initial datum $f^0$
satisfies
$$
\sup_{\tau \geq 0} \int_{\rr^3} \vert v \vert^4 \, f(t, v) \,
dv<\infty.
$$
\end{lemma}

\proof Without loss of generality we can assume that $f^0$, and
hence $f(t)$ for all $ \tau \geq 0$, has zero mean velocity and
unit temperature. We let
$$
m_4(t) = \int_{\rr^3} \vert v \vert^4 \, f(t, v) \, dv
$$
denote the fourth order moment of $f(t)$. Then, using the weak
formulation of the inelastic Boltzmann equation, we have:
\begin{equation}\label{eq-mom4}
\frac{d m_4(t)}{d t} =  \int_{\rr^3} \vert v \vert^4 \,
\tilde{Q}_e(f(t), f(t)) (v) \, dv
\end{equation}
that can be computed as in \cite{Bolley-Carrillo} by
\begin{align*}
\int_{\rr^3} \vert v \vert^4 \, \tilde{Q}_e(f, f) (v) \, dv = &-
<\zeta> \, \int_{\rr^3} \vert v \vert^4 \, f(v) \, dv   + <\mu_1>
\Big( \int_{\rr^3} \vert v \vert^2 \,  f(v) \, dv \Big)^2
\\ &+ <\mu_2> \iint_{\rr^3 \times \rr^3} (v \cdot w)^2 \, f(v) \,
f(w) \, dv \, dw
\end{align*}
where
$$
\mu_1 = \frac{1}{8} (\nu_1 + \nu_2 - \nu_3) \quad {\textrm{and}}
\quad \mu_2 = \frac{1}{4} (\nu_1 - \nu_2)
$$
with
$$
\nu_1 = (\epsilon^2 + \epsilon'^2)^2 -1 + \frac{4}{3} \epsilon^2
\epsilon'^2, \quad  \nu_2 = 2 \big[ \epsilon^2 + \epsilon'^2 -1 +
\frac{2}{3} \epsilon'^2 \big], \quad \nu_3 = 4 (\epsilon^2 - 1),
$$
and
$$
\zeta =  \frac{1}{3} ( 1 + 4 \, \epsilon - 7 \, \epsilon^2 + 4 \,
\epsilon^3 - 2 \, \epsilon^4) \qquad \mbox{with} \qquad \epsilon =
\frac{1-\et}{2} \qquad \mbox{and} \qquad \epsilon'=1-\epsilon .
$$
Now, \eqref{eq-mom4} reads
\begin{equation}\label{new}
\frac{d m_4(t)}{d t} =  -<\zeta> m_4(t) + m(t)
\end{equation}
where $m(t)$ is a combination of second order moments, which are bounded in
time since the kinetic energy is preserved by equation \eqref{him}. Moreover
one can check from the expression of $\zeta$ in terms of $e$ that $<\zeta> = 1-
\A(2,e,\eta)>0$. This ensures that $m_4(t)$ is bounded uniformly in time if
initially finite, and concludes the argument.
\endproof

\

The preceding result also shows the divergence of the fourth moment in case the
random variable $\eta$ and the restitution coefficient $e$ are chosen to
satisfy $\A(2,e,\eta)
> 1$ but $\A(\alpha,e,\eta) < 1$ for some $0<\alpha<2$.

\begin{corollary}
Let the restitution coefficient $e$ and the random variable $\eta$
be chosen so that $\A(2,e,\eta)>1$ but $\A(\alpha,e,\eta)<1$ for
some $0<\alpha<2$. Then, the unique isotropic steady state
$f_\infty$ in $\Ptwot$ of equation~\eqref{him} with zero mean
velocity and unit pressure tensor has unbounded fourth moment.
\end{corollary}

\proof With the notation of the previous subsection, the evolution of the
fourth moment for isotropic densities given in Lemma \ref{prop-mom4} ensures
that
$$
\frac{d m_4(t)}{d t} =  -<\zeta> m_4(t) + m(t),
$$
where $m(t)$, which is a combination of second order moments, is bounded from
below.  Recall that $<\zeta>=1-\A(2,e,\eta)<0$ to conclude.
\endproof


\bigskip

\noindent {\bf Acknowledgements:} JAC acknowledges the support from DGI-MEC
(Spain) FEDER-project MTM2005-08024 and 2005SGR00611. G.T. acknowledges the
support of the Italian MIUR project ``Kinetic and hydrodynamic equations of
complex collisional systems''. JAC and GT acknowledge partial support of the
Acc. Integ. program HI2006-0111. JAC acknowledges partial support of the Acc.
Integ. program HF2006-0198.


\end{document}